\documentclass{amsproc}

\usepackage[english]{babel}

\newtheorem{thm}{Theorem}[subsection]
\newtheorem{prop}[thm]{Proposition}
\newtheorem{lem}[thm]{Lemma}
\newtheorem{cor}[thm]{Corollary}
\theoremstyle{definition}
\newtheorem{dfn}[thm]{Definition}
\theoremstyle{remark}
\newtheorem{rmk}[thm]{Remark}
\newtheorem{example}[thm]{Example}

\numberwithin{equation}{subsection}

\newcommand{\C}{{\mathbb C}}

\newcommand{\R}{{\mathbb R}}
\newcommand{\Z}{{\mathbb Z}}

\newcommand{\CF}{\mathcal F}

\newcommand{\CL}{\mathcal L}

\newcommand{\CO}{\mathcal O}

\newcommand{\ka}{\mathfrak a}

\newcommand{\kg}{\mathfrak g}  
\newcommand{\kh}{\mathfrak h}
\newcommand{\kk}{\mathfrak k}  
\newcommand{\kl}{\mathfrak l}
\newcommand{\km}{\mathfrak m}
\newcommand{\kn}{\mathfrak n}
\newcommand{\ko}{\mathfrak o}
\newcommand{\kp}{\mathfrak p}  
\newcommand{\kq}{\mathfrak q}
\newcommand{\kt}{\mathfrak t}
\newcommand{\ku}{\mathfrak u}
\newcommand{\ks}{\mathfrak s}

\DeclareMathOperator{\ad}{ad}

\newcommand{\ga}{\alpha}
\newcommand{\gb}{\beta}
\newcommand{\gga}{\gamma}           \newcommand{\gG}{\Gamma}
           \newcommand{\gD}{\Delta}

\newcommand{\gth}{\theta}

\newcommand{\gs}{\sigma}

\newcommand{\go}{\omega}           \newcommand{\gO}{\Omega}

\begin{document}

\title{Weighted Vogan diagrams associated to real nilpotent orbits}
\author{Esther Galina}
\address{CIEM-FAMAF, Universidad Nacional de Córdoba, Ciudad Universitaria, 5000 Córdoba, Argentina}
\email{galina\@famaf.unc.edu.ar}
\thanks{The author was supported in part by CONICET PIP 6308, SECYT-UNC, FONCYT}
\begin{abstract}
We associate to each nilpotent orbit of a real semisimple Lie algebra $\kg_o$ a weighted Vogan diagram, that is a Dynkin diagram with an involution of the diagram, a subset of painted nodes and a weight for each node. Every nilpotent element of $\kg_o$ is noticed in some subalgebra of $\kg_o$. In this paper we characterize the weighted Vogan diagrams associated to  orbits of noticed nilpotent elements.
\end{abstract}
\maketitle

\section*{Introduction}

Let $G_o$ a real connected reductive Lie group and $\kg_o$ its Lie
algebra.  The  purpose of this work is to associate to each
nilpotent $G_o$-orbit of a simple real Lie algebra $\kg_o$ a
diagram. It is known that  the problem can really be reduced to
study semisimple groups of adjoint type.

In the case where the group and the algebra are complex, the
Jacobson and Morozov theorem  relates the orbit of a nilpotent
element $e$ with a triple $(h,e,f)$ that generates a subalgebra of
$\kg$ isomorphic to $\frak{sl}(2,\C)$. There is a parabolic
subalgebra associated to this triple that permits to attach a
weight to each node of the Dynkin diagram of $\kg$. The resulting
diagram is called a weighted Dynkin diagram associated to the
nilpotent orbit of $e$ (see \cite{C}). The  Bala and Carter
classification \cite{BC} (see also \cite{C}) of  nilpotent
G-orbits of a complex reductive Lie algebra $\kg$ establishes a
one-to-one correspondence between nilpotent orbits of $\kg$ and
conjugacy classes of pairs $(\km, \kq_\km)$, where $\km$ is a Levi
subalgebra of $\kg$ and $\kq_\km$ is a distinguished parabolic
subalgebra of the semisimple Lie algebra $[\km,\km]$. The
nilpotent orbit comes from the Richardson orbit of the connected
subgroup $Q_\km$ of $G$ with Lie algebra $\kq_\km$. The elements
of the nilpotent orbit result distinguished in $\km$. Therefore,
it is of great importance the classification of the distinguished
parabolic subalgebras, it is done using the  weighted Dynkin
diagrams (see \cite {BC} for this classification for complex
simple Lie algebras).

In the case where the group and the algebra are real, consider a
Cartan decomposition  of $\kg_o=\kk_o+ \kp_o$. By complexification
there is a decomposition $\kg=\kk +\kp$ corresponding to a Cartan
involution $\theta$. Denote by $\sigma$ the conjugation in $\kg$
with respect to $\kg_o$. In this setting Sekiguchi \cite {S}
proves a  one-to-one correspondence, conjectured by Kostant,
between $G_o$-orbits in $\kg_o$ and $K$-orbits of $\kp$ where $K$
is the connected subgroup of the adjoint group $G$ of $\kg$ with
Lie algebra $\kk$. Following the ideas of the Bala-Carter
classification and using the Kostant-Sekiguchi correspondence,
N\"oel \cite{N} gives a classification of $G_o$-nilpotent orbits
in $\kg_o$. Actually, the classification was previously known but
 based on other circle of ideas, see \cite{CMc} for a complete
version that includes all the cases, or \cite{Kaw} for a general
analysis. N\"oel proves that the orbits $K.e$ are in one-to-one
correspondence with the triples of the form $(\km,\kq_\km,\kn)$
where $e$ is a non-zero nilpotent element of $\kp$, $\km$ is a
minimal $\theta$-stable Levi subalgebra of $\kg$ containing $e$,
$\kq_\km$ is a $\theta$-stable parabolic subalgebra of $[\km,\km]$
and $\kn$ is a certain $(L\cap K)$-prehomogeneous space of
$\kq_\km \cap \kp$ containing $e$, where $L$ is the Levi subgroup
of the corresponding parabolic subgroup of $G$ with Lie algebra
$\kq_\km $. In doing so N\"oel defines the {\it noticed} nilpotent
elements, it results that every nilpotent element $e$ is noticed
in the minimal $\theta$-stable Levi subalgebra of $\kg$ containing
$e$. In analogy with distinguished nilpotent orbits corresponding
to distinguished parabolic subalgebras, in our situation noticed
nilpotent orbits of $\kp$ are in correspondence with noticed
triples $(\kg, \kq, \kn)$  as above. This explain the importance
of the classification of noticed nilpotent $K$-orbits in $\kp$.

In this paper we attach to each nilpotent $K$-orbit a weighted
Vogan diagram. It consists in a Vogan diagram (see \cite{K}) with
weights attached to the nodes. The values of the weights are in
the set $\{0,1,2\}$. Forgetting the painted nodes and the
involutive automorphism of the diagram, it is a weighted Dynkin
diagram. Moreover, from an abstract weighted Vogan diagram one can
re-obtain the real algebra $\kg_o$  and a triple $(\kg, \kq,
\kn)$. That is, there is an assignment from nilpotent
$K$-orbits in $\kp$ to equivalence classes of weighted Vogan
diagrams.

We intent to determine all the weighted Vogan diagrams associated
to noticed nilpotent $K$-orbits  in $\kp$. For classical real Lie
algebras there is a parameterization of nilpotent $G_o$-orbits of
$\kg_o$ by signed Young diagrams \cite{SS}, \cite{BCu},
\cite{CMc}, but not for exceptional ones. N\"oel determines all
the noticed orbits using this parameterization case by case
\cite{N}. For exceptional real Lie algebras he uses the
Djokovi\'c's tables  of the reductive centralizer of real
nilpotent orbits \cite{D1}, \cite{D2}. The classification of
noticed nilpotent $K$-orbits in $\kp$ by weighted Vogan diagrams
will give a classification of all noticed nilpotent orbits of a
simple real Lie algebra. In this paper we gives a characterization
of the weighted Vogan diagrams of noticed nilpotent orbits. This
characterization  will permit us to give a classification of all
noticed nilpotent orbits of a  simple real Lie algebra. It will be
the content of a future paper.

Kawanaka also gives a parameterization of nilpotent orbits of a simple real Lie algebra using weighted Dynkin diagrams (see \cite {Kaw}). It seems to be compatible with our description, but here we can explicitly reconstruct the real nilpotent orbit from a weighted Vogan diagram.

I want to thank Sebasti\'an Simondi for the helpful conversations
we have  had that developed into the Definition 3.2.1, and to the
referee of this work for the valuable observations.

\section{Real nilpotent orbits} \label{data}

\subsection{}\label{thetasigma} Let $G_o$ be a real semisimple Lie group of adjoint type with Lie algebra $\kg_o$. Let $\theta$ a Cartan involution of $\kg_o$. It gives place to a Cartan decomposition in eighenspaces of $\kg_o = \kk_o \oplus \kp_o$ where $\kk_o$ is the subalgebra of $\theta$-fixed points. Denote the complexification of a space by the same letter but without the subscript. Extend the Cartan involution $\theta$ to $\kg$ linearly.  Let $\gs$ be the conjugation in $\kg$ with respect to $\kg_o$.

If $G$ is the adjoint group of $\kg$, denote by $K$ the connected subgroup with Lie algebra $\kk$.

\subsection{}\label{JM} The Jacobson-Morosov theorem associates to each nilpotent element $e\in \kg$ a triple (or JM-triple) $(h,e,f)$ such that $[h,e]=2e$, $[h,f]=-2f$ and $[e,f]=h$ (see \cite{C}). That is, the subalgebra of $\kg$ generated by the triple is isomorphic to $\mathfrak{sl}(2,\C)$. Moreover, there is a one-to-one correspondence between non-zero nilpotent $G$-orbits in $\kg$ and $G$-conjugacy classes of the Lie subalgebras generated by a triple of $\kg$ of this kind.

\subsection{} Following the work of Kostant-Rallis \cite{KR} if $e\in \kp$ there exists a JM-triple $(h,e,f)$ with $h\in K$ and $f\in \kp$ called a {\it normal} triple. In this case the one-to-one correspondence is between non-zero nilpotent $K$-orbits in $\kp$ and $K$-conjugacy classes of the Lie subalgebras generated by a normal triple of $\kg$.

\subsection{}\label{KStriple} On the other hand, Sekiguchi obtained in \cite{S} that each $K$-conjugacy class of  the Lie subalgebras generated by a normal triple of $\kg$ contains a subalgebra generated by a normal triple such that $f=\gs(e)$ and $h\in i\kk_o$. Following N\"oel, we called it a {\it KS-triple} (KS comes from Kostant and Sekiguchi). That is, each nilpotent element of $\kp$ is $K$-conjugate to an element $e$ of a KS-triple $(h,e,f)$.

Moreover, he proves that for each nilpotent $G_o$-orbit $\CO$ of $\kg_o$ there is a real JM-triple $(h_o,e_o,f_o)$ that generates a subalgebra isomorphic to $\mathfrak{sl}(2,\R) $ in $\kg_o$ such that  $e_o\in\CO$ and $\theta(h_o)=-h_o$ and $\theta(e_o)=-f_o$. In this setting, the triple $(h,e,f)$ in $\kg$
\begin{align*}
h=c(h_o)&= i(e_o-f_o) \in i\kk_o  \\
e=c(e_o)&= \frac 12(h_o-i(e_o+f_o)) \\
f=c(f_o)&= \frac 12(h_o+i(e_o+f_o)),
\end{align*}
given by the Cayley transform $c$, is a KS-triple.

The application $G_o.e_o \to K.e$ produce a one-to-one correspondence between nilpotent $G_o$-orbits in $\kg_o$ and nilpotent $K$-orbits in $\kp$.

\section{Nilpotent orbits of real symmetric pairs}  \label{sec:2}

Continue with the notation of previous section. N\"oel in \cite{N} gives a parameterization of nilpotent $K$-orbits in $\kp$ following the philosophy of Bala-Carter classification for nilpotent $G$-orbits in $\kg$. According to the Kostant-Sekiguchi correspondence of the previous section, this results a classification of nilpotent $G_o$-orbits in $\kg_o$.

\subsection{}\label{g(i)} It is known that each JM-triple $(h,e,f)$ in $\kg$, determines a graduation
$$
\kg= \oplus_{i\in \Z} \ \kg^{(j)}
$$
where  $\kg^{(j)}= \{x\in \kg : [h,x]= j\, x\}$ (see \cite {BC} for more details). Evidently $h\in \kg^{(0)}$,  $e\in \kg^{(2)}$ y $f\in \kg^{(-2)}$.
These eigenspaces have the following property
$$
[\kg^{(i)} ,\kg^{(j)} ]\subset \kg^{(i+j)}
$$
for all $i,j\in \Z$.

Denote by $\kl = \kg^{(0)}$, $\ku= \oplus_{j>0} \ \kg^{(j)}$ and
$\bar\ku= \oplus_{j<0} \ \kg^{(j)}$. These subspaces are subalgebras of $\kg$. The direct sum $\kq = \kl \oplus \ku$ is a parabolic subalgebra that contains $e$. In fact, as $\ad h$ is semisimple in $\kg$, there is a Cartan subalgebra $\kh$ that contains  $h$.  Let $\Psi=\{\ga_1,\dots,\ga_n\}$ be a set of simple roots of the root system $\gD(\kg,\kh)$ determined by the condition $\gb(h)\geq 0$ for all $\gb \in \Psi$. Notice that there are several choices of $\Psi$ with this condition. In fact, if  $W_o$ is the subgroup of the Weyl group of $\kg$ generated by the roots $\gb$ that satisfies $\gb(h)=0$,  $w\Psi$ is another set of simple roots with the same condition for any $w\in W_o$. Therefore,
$$
\kg^{(j)}= \sum_{\ga(h)=j}\ \kg_\ga
$$
where $\kg_\ga$ are the corresponding rootspaces. So,
$$
\kq = \kl \oplus \ku = \kg^{(0)} \oplus \sum_{j>0}\ \kg^{(j)}
=\kh \oplus \sum_{\ga(h)\geq 0}\ \kg_\ga
$$
The subalgebra $\kq=\kq_{(h,e,f)}$ is called the {\it parabolic subalgebra associated to the JM-triple $(h,e,f)$}.

\begin{rmk} \label{rmk:thetastable}
If $(h,e,f)$ is a KS-triple, choose the Cartan subalgebra $\kh$ in
the following way. As $h \in\kk$, take $\kt$ a maximal abelian
subspace of $\kk$ that contains $h$ and  a  maximal abelian
subspace  $\ka$ in $\kp$ that conmutes with $\kt$. So $\kh= \kt
\oplus \ka$ is  $\gth$-stable and is a maximal compact Cartan
subalgebra of $\kg$ that contains $h$. Therefore, the parabolic
subalgebra $\kq$,  the root system $\gD(\kg,\kh)$ and the set of
simple roots $\Psi$ result $\gth$-stable (see \cite{N} for more
details).
\end{rmk}

Denote by $Q$ and $L$ the analytic subgroups of $G$ with Lie algebras $\kq$ and $\kl$ respectively.

\begin{prop} \label{prop:g1}
Let $(h,e,f)$ be a KS-triple, then

$(i)$ $(L\cap K).e$ is dense in $\kg^{(2)}\cap \kp$;

$(ii)$ $(Q\cap K).e $ is dense in $\oplus_{i\geq 2}\ \kg^{(i)}\cap \kp$;

$(iii)$ $\dim \kg^{(1)}\cap \kk =  \dim \kg^{(1)}\cap \kp$.
\end{prop}

The first item of the proposition was proved by Kostant and Rallis (\cite {KR}, in proof of Lemma 4). The two other are results of N\"oel \cite {N}, \cite {N2}.

%

\subsection{} \label{wDd}
The correspondence between nilpotent $G$-orbits of $\kg$ and
$G$-conju\-ga\-cy classes of the Lie algebras generated by a
JM-triple of $\kg$ permits to associate to each nilpotent
$G$-orbit a {\it weighted Dynkin diagram}. It consists in a pair
$(D,\go)$ where $D$ is the Dynkin diagram of $\kg$ and  $\go$ is a
set of weights attached to the nodes of the diagram. If $(h,e,f)$
is a JM-triple corresponding to a nilpotent $G$-orbit, $\go$ is
defined by $\go_i= \ga_i(h)$ where $\Psi =\{\ga_1,\dots,\ga_n\}$
is the set of simple roots of $\gD(\kg,\kh)$ defined in
\ref{g(i)}. Note that two weighted Dynkin diagrams associated to a
pair of JM-triples are equal if and only if the triples are in the
same $G$-conjugacy class. This means that the weighted Dynkin
diagram only depends on the nilpotent $G$-orbit (see \cite {BC} or
\cite {C} for more details).

\subsection{} We need some definitions to explicitly enunciate the results of N\"oel. In
Bala-Carter results {\it distinguished} nilpotent elements play an
important role. They are  defined as the nilpotent elements whose
centralizers do not contain any semisimple element.  Or
equivalently, $e$ is distinguished if the minimal Levi subalgebra
of $\kg$ that contains it is $\kg$ itself. Classification of
weighted Dynkin diagrams of distinguished nilpotent $G$-orbits
gives a parametrization of nilpotent $G$-orbits because each
nilpotent element in $\kg$ is distinguished in the minimal Levi
subalgebra that contains it. For symmetric pairs this role is
played by noticed nilpotent elements.

\begin{dfn} \label{dfn: noticed} {\rm (N\"oel \cite {N})}
A nilpotent element $e\in \kp$ is {\it noticed} if $\kg$ is the
minimal $\gth$-stable  Levi subalgebra that contains $e$.

\noindent We will say that a KS-triple $(h,e,f)$  is {\it noticed} if $e$ is noticed, as well as the nilpotent $K$-orbit $K.e$ and the real nilpotent $G_o$-orbit associated to $K.e$ by the Kostant-Sekiguchi correspondence.
\end{dfn}

\rmk \label{rmk:distinguished}
{Every distinguished nilpotent element in $\kp$ is noticed but the converse is not true.}

\begin{rmk} \label{rmk:notLevi} {Every nilpotent element in $\kp$ is noticed in the minimal $\gth$-stable Levi subalgebra that contains it.}
\end{rmk}

As we consider KS-triples we are in  the situation of Remark \ref{rmk:thetastable}.

\begin{prop} \label{prop:noticed}{\rm (N\"oel \cite {N})}
The following statements are equivalent:

$(i)$ $e\in \kp$ is noticed;

$(ii)$ the centralizer $\kk^{(h,e,f)}$ of the noticed triple
$(h,e,f)$ in $\kk$ is 0.

$(iii)$ $\dim \kg^{(0)}\cap \kk =  \dim \kg^{(2)}\cap \kp$.
\end{prop}

\rmk {Every distinguished nilpotent element is even, that is $\kg^{(i)}=0$ for all integer $i$ odd, but this is not true for noticed nilpotent elements.

\subsection{}\label{g1}
We will express in other terms some of previous results, not given
so explicitly in \cite{N}.  Fix a KS-triple $(h,e,f)$ and continue
with the same notation and considerations of \ref{g(i)}. Each
$\ga$ in the set of positive roots $\gD^+$ associated to the set
of simple roots $\Psi=\{\ga_1,\dots,\ga_n\}$ is of the form
$\ga=\sum_{i=1}^n n_i \ga_i$ for certain non negative integers
$n_i$. For each $\ga \in \gD^+$ define its weight $ \go_\ga=
\ga(h)$ and the number
$$
l_\ga = \sum_{{\ga_i}\in \Psi} n_i \qquad \qquad m_\ga =
\sum_{\kg_{\ga_i}\in \kp} n_i
$$
We will call them the {\it lenght} and {\it non-compact lenght} of
$\ga$ respectively. Denote by
$$\begin{aligned}
M_\kk^{(j)} &= \{\ga \in \gD^+ : \go_\ga=j,\ \gth\ga=\ga,\ m_\ga \text{ even},\ \ga\neq \gga+\gth\gga \text{ for some } \gga \in \gD^+\}\\
M_\kp^{(j)} &= \{\ga \in \gD^+ : \go_\ga=j,\ \gth\ga=\ga, \ m_\ga \text{ odd}\} \cup \{\gga+\gth\gga \in \gD^+ : \gga \in \gD^+,\ 2\go_\gga=j\}\\
C^{(j)} &= \{ \{\ga,\ \gth\ga\} : \ga \in \gD^+,\ \gth\ga\neq \ga,\ \go_\ga=j\}\\
\end{aligned}
$$

\begin{rmk}
Note that the only situation of $\kg$ simple where the second set
of  $M_\kp^{(j)}$ is  not zero is  in the case $\kg=\ks\kl(2m,\C)$
and $\gth$ does not fix any simple root. The reason is that the
only Dynkin diagram with an automorphism of order two such that
$\gga+\gth\gga$ is a root for some $\gga \in \gD^+$ is the case
mentioned above. The weight of the root $\gga+\gth\gga$ is
$2\go_{\gga}$, so $j$ is even in this case.
\end{rmk}

\begin{prop} \label{lem:g(j)}
The spaces $\kg^{(j)}\cap \kk$ and $\kg^{(j)}\cap \kp$ with $j\geq 0$ can be described as
$$\begin{aligned}
\kg^{(0)}\cap \kk &= \kt + \sum_{\ga\in \pm M_\kk^{(0)}} \kg_{\ga}
\oplus
\sum_{\{\ga, \gth \ga\}\in \pm C^{(0)}} \C(X_{\ga}+\gth X_{\ga})\\
\kg^{(j)}\cap \kk &= \sum_{\ga\in M_\kk^{(j)}} \kg_\ga \oplus
\sum_{\{\ga, \gth \ga\}\in C^{(j)}} \C(X_\ga+\gth X_\ga) \qquad j>0\\
\kg^{(0)}\cap \kp &= \ka + \sum_{\ga\in \pm M_\kp^{(0)}} \kg_\ga \
\oplus
\sum_{\{\ga, \gth \ga\}\in \pm C^{(0)}} \C(X_\ga-\gth X_\ga) \\
\kg^{(j)}\cap \kp &= \sum_{\ga\in M_\kp^{(j)}} \kg_\ga \ \oplus
\sum_{\{\ga, \gth \ga\}\in C^{(j)}} \C(X_\ga-\gth X_\ga) \qquad j>0\\
\end{aligned}
$$
\end{prop}

\begin{proof}
Without lost of generality we can consider $\kg$ simple.

Let $\kg_\ga$ a root space such that $\go_\ga=j>0$. If $\gth\ga\neq \ga$, then $\kg_\ga\oplus \gth\kg_\ga= \kg_\ga\oplus \kg_{\gth\ga}$ is $\gth$-stable and is generated by $\{X_\ga+\gth X_\ga,X_\ga-\gth X_\ga\}$ for some $0\neq X_\ga\in \kg_\ga$. The first generates $(\kg_\ga\oplus \kg_{\gth\ga})\cap \kk$ and the second $(\kg_\ga\oplus \kg_{\gth\ga})\cap \kp$. So, we have obtained the $\ga$-component in the last sumands of the l.h.s. of the equalities.

 If $\gth\ga = \ga$, let's see that $\ga=\sum_{i=1}^n n_i \ga_i$ is compact or non compact depending on the parity of $m_\ga$. Let's prove it by induction on the length $l_\ga=\sum_{i=1}^n n_i$ of $\ga$.

For $l_\ga=1$ it is obvious. Consider $l_\ga >1$. Hence there
exists a root  $\gb$ fixed by $\gth$ such that $\ga = \gb + \ga_s$
for some $s$ if $\gth \ga_s=\ga_s$, or $\ga = \gb + \ga_s + \gth
\ga_s$ if $\gth \ga_s\neq \ga_s$. In the last case, it is possible
that $\gb=0$.

In the first case, $ \kg_\ga \subset [\kg_{\gb},\kg_{\ga_s}]. $
So, if $\kg_{\ga_s} \in \kp$, $m_\ga = m_\gb + 1$. By inductive
hypotesis, $m_\gb$ odd implies that $\kg_{\gb} \in \kp$. Then,
$m_\ga $ is even and $\kg_{\ga} \in \kk$. The case $m_\gb$ even
implies that $\kg_{\gb} \in \kk$. Then, $m_\ga $ is odd and
$\kg_{\ga} \in \kp$. If $\kg_{\ga_s} \in \kk$, $m_\ga = m_\gb$.
So, $\kg_{\ga} \in \kk$ if and only if $\kg_{\gb}\in \kk$.

In the second case, suppose $\gb\neq 0$. Hence, $\kg_{\ga_s} + \gth\kg_{\ga_s}$ is not a root, according with the possible automorphism of Dynkin diagrams of order two. Therefore,
$$
\kg_\ga \subset
[\kg_{\ga_s},[\kg_{\gb},\kg_{\gth\ga_s}]]
$$
and $m_\ga = m_\gb$. Consider $X_{\ga_s} \in \kg_{\ga_s}$ and
$X_{\gb} \in \kg_{\gb}$. As $\kg_{\gth\ga_s} = \gth \kg_{\ga_s}$ ,
we can analyse
$$
\gth[X_{\ga_s},[X_{\gb},\gth X_{\ga_s}]]=
[\gth X_{\ga_s},[\gth X_{\gb},X_{\ga_s}]]=
[[\gth X_{\ga_s},\gth X_{\gb}],X_{\ga_s}] + [\gth X_{\gb},[\gth  X_{\ga_s},X_{\ga_s}]].
$$
But, $[\gth  X_{\ga_s},X_{\ga_s}] =0$. Hence,
$$
\gth[X_{\ga_s},[X_{\gb},\gth X_{\ga_s}]]=
[[\gth X_{\ga_s},\gth X_{\gb}],X_{\ga_s}]=
[X_{\ga_s},[\gth X_{\gb},\gth X_{\ga_s}]]
$$
Then, $\kg_{\ga} \in \kk$ if and only if $\kg_{\gb} \in \kk$.

If $\gb=0$, it means that ${\ga_s} + \gth{\ga_s}$ is a root. Then, there exits $X_{\ga_s}$ such that $0 \neq [X_{\ga_s},\gth X_{\ga_s}] \in \kg_{{\ga_s} + \gth{\ga_s}}$ and
$$
\gth[X_{\ga_s},\gth X_{\ga_s}]= [\gth X_{\ga_s}, X_{\ga_s}]=-[X_{\ga_s},\gth X_{\ga_s}]
$$
This says that $\kg_\ga \in \kp$. So, we have obtained the
$\ga$-component in the first sumands of the  l.h.s. of the
equalities. Therefore, the equalities hold for $j>0$. If $j=0$, we
also have to consider the subspaces $\kt$ and $\ka$ respectively,
and the sets of negative roots $-M_\kk^{(0)}$, $-M_\kp^{(0)}$ and
$-C^{(0)}$ that have weight zero too.
\end{proof}

Given a set $U$, denote by $|U|$ the cardinality of $U$.

\begin{cor}\label{cor:M1}
Let $(h,e,f)$ be a KS-triple of $\kg$, then the sets $M_\kk^{(1)}$
and $M_\kp^{(1)}$ associated to it have the same cardinality.
\end{cor}

\begin{proof}
This follows immediatly from Proposition \ref{prop:g1} $(iii)$
because
$$
0=\dim \kg^{(1)}\cap \kk - \dim \kg^{(1)}\cap \kp = |M_\kk^{(1)}|
+ |C^{(1)}|- (|M_\kp^{(1)}|+ |C^{(1)}|)=|M_\kk^{(1)}| -
|M_\kp^{(1)}|
$$
\end{proof}

\begin{cor}\label{cor:d0k=d2p}
Let $(h,e,f)$ be a KS-triple of $\kg$, it is noticed if and only
if $\dim\kt + 2|M_\kk^{(0)}| + 2|C^{(0)}|= |M_\kp^{(2)}| +
|C^{(2)}|$.
\end{cor}

\begin{proof}
This follows inmediatly from Proposition \ref{prop:noticed}
$(iii)$ and the previous proposition.
\end{proof}

\subsection{} Let $B$ be the Killing form of $\kg$. Define $B'$ by
$B'(x,y)= -B(x,\gth\gs(y ))$ for all $x,y\in\kg$. It results a definite positive hermitian form.

Let $\kq = \kl \oplus \ku$ be a $\gth$-stable parabolic subalgebra of $\kg$, that is $\kl$ and $\ku$ are. Let $\mu$ be the subspace of $\ku\cap\kp$ such that the decomposition
$$
\ku\cap\kp = \mu \oplus [\ku\cap\kk, [\ku\cap\kk,\ku\cap\kp ]]
$$
is orthogonal with respect to $B'$. Let $L$ be the analytic subgroup of $G$ with Lie algebra $\kl$. Let $\eta$ be an $(L\cap K)$-module of $\mu$ and
$\hat\eta = \eta \oplus [\kl\cap \kp,\eta]$; $\hat\eta$ is a $\gth$-stable subspace of $\kg$.

Define $\CL_\kg$ the set of pairs $(\kq,\eta)$ where $\kq=\kl \oplus \ku$ is a $\gth$-stable parabolic subalgebra of $\kg$ and $\eta$ as above such that
\begin{enumerate}
\item $\eta$ has a dense $(L\cap K)$-orbit;
\item $\hat\eta$ has a dense $L$-orbit;
\item $\dim \kl\cap \kk =  \dim \eta$;
\item $\hat\eta$ is orthogonal to $[\ku,[\ku,\ku]]$;
\item $\hat\eta$ is orthogonal to $[\ku,\hat\eta]$;
\item $[\ku\cap\kk,\ku\cap\kp ]\subset [\kq\cap\kk,\eta]$.
\end{enumerate}





\begin{rmk}
If $(h,e,f)$ is a KS-triple, then the decomposition (2.1) of $\kg$ in eigen\-spa\-ces of $\ad h$ of $\kg$ is orthogonal with respect to $B'$ and the parabolic subalgebra associated to it is $\gth$-stable.
\end{rmk}

Having in mind \ref{KStriple} and previous definitions, we can state an important correspondence.

\begin{thm}\label{thm:noticedcorresp} {\rm (N\"oel \cite {N})} {There is a one-to-one correspondence between $K$-conjugacy classes of Lie subalgebras generated by noticed KS-triples and  $K$-conjugacy classes of pairs $(\kq,\eta) \in \CL_\kg$. }
\end{thm}

The map
$(h,e,f) \to (\kq_{(h,e,f)},\kg^{(2)}\cap \kp)$
sends $K$-conjugacy classes of noticed KS-triples into  $K$-conjugacy classes of pairs $(\kq, \eta) \in \CL_\kg$. It inverse is defined choosing a nilpotent element of the dense $(L\cap K)$-orbit of the space $\eta$ and considering a KS-triple associated to it.

The main theorem of this work of N\"oel is a consequence of Theorem \ref{thm:noticedcorresp} considering Remark \ref{rmk:notLevi}.

\begin{thm} \label{thm:Noelcorresp} {\rm (N\"oel \cite {N})}
There is a one-to-one correspondence between $K$-orbits of
nilpotent elements of $\kp$ and $K$-conjugation clases of pairs
$(\kq_\km,\eta_\km) \in \CL_\km$ with $\km$ running over all
$\gth$-stable Levi subalgebras of $\kg$.
\end{thm}

\section{Abstract weighted Vogan diagrams}

 The purpose of this section is to define abstract weighted Vogan diagram. In the next one we will relate them with nilpotent $K$-orbits in $\kp$.

\subsection{} Here we combine the notions of Vogan diagrams, used in the classification of real simple Lie algebras \cite {K}, and weighted Dynkin diagram, used to classify nilpotent $G$-orbits \cite {BC}, \cite {C}.

\begin{dfn} An  {\it abstract Vogan diagram} is a diagram with data $(D,\theta,J)$ where $D$ is a Dynkin diagram of $n$ nodes, $\theta$ is an automophism of the diagram $D$ of order at most 2 and $J$ is a subset of $\theta$-invariant nodes of $D$.

\noindent An {\it abstract weighted Vogan diagram} $(D,\theta,J, \go)$ consists on an abstract Vogan diagram  $(D,\theta,J)$ with a set of weights $\go= (\go_1,\dots,\go_n)$ associated to the nodes that satisfy $\go_i = \go_{\theta(i)}$  and $\go_i \in \{0,1,2\}$ for all $i=1,\dots,n$.
\end{dfn}

\begin{rmk} A weighted Vogan diagram $\gG=(D,\theta,J,\go)$  gives place to the weighted Dynkin diagram $(D,\go)$ forgetting the automorphism and the painted nodes \cite {C}.
\end{rmk}

In order to draw the diagram, if $\theta$  have orbits of 2 elements, the nodes in the same orbit are connected by a doublearrow. The nodes in the set $J$ are painted and each weight is written above the corresponding node.
For example, the diagram of Figure \ref{example-fig} corresponds to the data
$D=D_6$, $\theta$ the unique non trivial automorphism of $D_6$ that fix the first four nodes, $J=\{1,4\}$ and $\go=\{2,0,0,0,1,1\}$.

\begin{figure}[h]\label{fig:D6}
\begin{center}
\setlength{\unitlength}{10pt}
\begin{picture}(12,7)(-6,-3)
  \thicklines
  \put(-6,0){\circle*{.6}}
  \put(-3,0){\circle{.6}}
  \put(0,0){\circle{.6}}
  \put(3,0){\circle*{.6}}
  \put(5,1.5){\circle{.6}}
  \put(5,-1.5){\circle{.6}}
  \put(-5.7,0){\line(3,0){2.5}}
  \put(-2.7,0){\line(3,0){2.5}}
  \put(0.3,0){\line(3,0){2.5}}
  \put(3,0){\line(5,4){1.7}}
  \put(3,0){\line(5,-4){1.7}}

  \thinlines
  \put(5,1.1){\vector(0,-1){2.2}}
  \put(5,-1.1){\vector(0,1){2.2}}

\put(-6,1){\small{2}}
  \put(-3,1){\small{0}}
  \put(0,1){\small{0}}
  \put(3,1){\small{0}}
  \put(5.5,2){\small{1}}
  \put(5.5,-1){\small{1}}

\end{picture}
\end{center}
\caption{}
\label{example-fig}
\end{figure}

\begin{rmk}
An automorphism of a Dynkin diagram of order two is unique up to an exterior automorphism of the diagram. More explicitly, if the diagram is connected, it is unique except for $D_4$.
\end{rmk}

\subsection{} \label{equivalent}
Vogan proved that to each abstract Vogan diagram $(D,\theta,J)$  one can associate a 4-tuple $(\kg_o,\gth,\kh_o,\gD^+)$ of a real Lie algebra $\kg_o$, a Cartan involution $\gth$ of $\kg$, a real $\theta$-stable maximally compact Cartan subalgebra $\kh_o= \kt_o \oplus \ka_o$ of $\kg_o$ and a positive root system $\gD^+(\kg,\kh)$ that takes $i\kt_o$ before $\ka_o$ (see \cite {K}, Theorem 6.88). This permits the classification of all simple real Lie algebras, but it is possible that two different abstract Vogan diagrams give place to the same simple real Lie algebra.

Given an abstract Vogan diagram $\gG=(D,\theta,J,\go)$ we will say that the 4-tuple $(\kg_o,\gth,\kh_o,\gD^+)$ is {\it associated to $\gG$} if it is the associated to $(D,\theta,J)$.

\begin{dfn} \label{dfn:equivalent}
An abstract weighted Vogan diagrams $\gG=(D,\theta,J,\go)$ is  {\it equivalent} to a second one if one can pass from $\gG$ to the other in  finite operations of the type:

(A) given $j\in J$ with $\go_j=0$, the resulting weighted Vogan diagram is $\gG'=(D,\theta,J',\go)$ where
$$
 J'=\{i\in J : i\text{ is not adjacent to } j\} \cup \{i \not \in J: i \text{ adjacent to $j$} \}
$$
Except in the cases:

${\bf B_n:}\ j=n$ $ J'= J$,

${\bf C_n:}\ j=n-1$,
$J'= \{i\in J : i\neq n-2\} \cup \{n-2 \text{ if } n-2 \not \in J \}$,

${\bf F_4:}\ j=2$,
$J'= \{i\in J : i\neq 1\} \cup \{1 \text{ if } 1 \not \in J \}$.
\end{dfn}

Given $j\in J$ with $\go_j=0$, operation (A) is nothing more than
change the  colors of the nodes adyacent to $j$, except for long
neighbors of short roots in types {\bf B, C} and {\bf F}. For
example, the two weighted Vogan diagrams of Figure \ref{equiv-fig}
are equivalent, one is obtained from the other applying operation
(A) on the second node.

\begin{figure}[h]\label{fig:equivalent}
\begin{center}
\setlength{\unitlength}{10pt}
\begin{picture}(12,7)(-6,-3)
  \thicklines
  \put(-9,0){\circle*{.6}}
  \put(-6,0){\circle*{.6}}
  \put(-3,0){\circle{.6}}
  \put(-8.7,0){\line(3,0){2.5}}
  \put(-5.7,-0.2){\line(3,0){2.5}}
  \put(-5.7,0.2){\line(3,0){2.5}}
  \put(-4.5,-0.2){$\rangle$}

  \put(-9,1){\small{1}}
  \put(-6,1){\small{0}}
  \put(-3,1){\small{1}}
 \thicklines
  \put(3,0){\circle{.6}}
  \put(6,0){\circle*{.6}}
  \put(9,0){\circle*{.6}}
  \put(3.3,0){\line(3,0){2.5}}
  \put(6.3,-0.2){\line(3,0){2.5}}
  \put(6.3,0.2){\line(3,0){2.5}}
  \put(7.5,-0.2){$\rangle$}

  \put(3,1){\small{1}}
  \put(6,1){\small{0}}
  \put(9,1){\small{1}}
\end{picture}
\end{center}
\caption{}
\label{equiv-fig}
\end{figure}

\subsection{}\label{abnoticed}
We will define the notion of noticed abstract weighted Vogan
diagram to be used later.

Given an abstract weighted Vogan diagram $\gG=(D,\theta,J,\go)$
denote by $N^\gth$ the number of nodes of $D$ fixed by $\gth$ and
by $N^{\gth}_2$ the number of $\gth$-orbits with two elements.
Consider $(\kg_o,\gth,\kh_o,\gD_\gG^+)$, the 4-tuple associated to
$\gG$. Denote by $\ga_j$ the simple root of $\gD_\gG^+$
corresponding to the node $j$. Hence, every root $\ga \in
\gD_\gG^+$ has the decomposition: $\ga=\sum_{i=1}^n n_i \ga_i$ for
certain non negative integers $n_i$. Its {\it weight} and its {\it
painted lenght} are respectively given by
$$
\go_\ga= \sum_{i=1}^n n_i w_i \qquad \qquad  p_\ga = \sum_{i\in J}
n_i
$$
Denote by
$$\begin{aligned}
P_{np}^{(j)}=P_{np}^{(j)}(\gG) &= \{\ga \in \gD_\gG^+ : \go_\ga=j,\ \gth\ga=\ga,\ p_\ga \text{ even},\ \ga\neq \gga+\gth\gga \text{ for some } \gga \in \gD_\gG^+\}\\
P_p^{(j)} =P_p^{(j)}(\gG)&= \{\ga \in \gD_\gG^+ : \go_\ga=j,\ \gth\ga=\ga, \ p_\ga \text{ odd}\} \cup \{\gga+\gth\gga \in \gD_\gG^+ : \gga \in \gD_\gG^+,\ 2\go_\gga=j\}\\
K^{(j)} =K^{(j)} (\gG)&= \{ \{\ga,\ \gth\ga\} : \ga \in \gD_\gG^+,\ \gth\ga\neq \ga,\ \go_\ga=j\}\\
\end{aligned}
$$

\begin{dfn} \label{dfn:equivalent}
An abstract weighted Vogan diagrams $\gG=(D,\theta,J,\go)$  is
{\it noticed} if the following equality holds for the
corresponding subsets of $\gD_\gG^+$,
$$
N^\gth +N^\gth_2+ 2|P_{np}^{(0)}| + 2|K^{(0)}|= |P_p^{(2)}| +
|K^{(2)}|
$$
\end{dfn}

\section {Real nilpotent orbits and weighted Vogan diagrams}

Let $K$ and $\kp$ corresponding to a symmetric pair as in Section 1. We will attach to each nilpotent $K$-orbit of $\kp$ a weighted Vogan diagram. The main result is a correspondence between classes of abstract weighted Vogan diagrams and real nilpotent orbits.

\subsection{} \label{wVdKStriple}
Let $(h,e,f)$ be a KS-triple of $\kg$. We will associate to it an
abstract weighted  Vogan diagram.  To this triple we can attach a
weighted Dynkin diagram following \cite {BC}, \cite {C}. It
consists in the Dynkin diagram of $\kg$ and weights in each node
defined by $\go_i= \ga_i(h)$  according to \ref{g(i)}.  On the
other hand, by  Remark \ref{rmk:thetastable} the Cartan involution
$\gth$ of $\kg$ given in \ref{thetasigma} provides an automorphism
of weighted Dynkin diagrams where $\theta(i)$ is the node
corresponding to $\gth\ga_i$. As $h\in i\kk$,  we have that
$\go_{\gth(i)}=\gth\ga_i(h)= \ga_i(\gth h)= \ga_i(h)$. That is,
$\go_{\gth(i)}= \go_i$ for all  $i$.  Observe that the weights do
not change if one replace the triple by a $K$-conjugate triple.
The KS-triple $(h,e,f)$ remains being a KS-triple by
$K_o$-conjugation, $K$-conjugation preserves normality. This means
that $\go$ only depends on the $K$-orbit.

To obtain an abstract weighted Vogan diagram  it remains to define the set $J$ of painted nodes. Define $J$ the set of nodes fixed by the automorphism $\gth$ that correspond to non compact roots, that is those roots $\ga$ such that $\kg_\ga \in \kp$. This fact is also invariant by $K$-conjugation.

Note that in such assignment there is a choice of the set of
simple roots  $\Psi=\{\ga_1,\dots, \ga_n\}$ and there are $|W_o|$
posibilities of this choice, where $W_o$ is the subset of the Weyl
group generated by the set of roots $\ga$ such that $\ga(h)=0$. On
the other hand, observe that two weighted Vogan diagrams $\gG_1$
and $\gG_2$ obtained  from a KS-triple $(h,e,f)$, like before,
associated to different sets of simple roots $\Psi_1$ and $\Psi_2=
w \Psi_1$ for some $w\in W_o$, are equivalent. In fact, $w$ is a
composition of finite reflections $s_{\ga_j}$ with $\ga_j$ a
simple root of $\Psi_1$ in $W_o$. The set $s_{\ga_j} \Psi_1$ gives
rise to the same  weighted Vogan diagram if $\ga_j$ is compact or
complex, or to one that can be obtained from $\gG_1$ applying
operation (A) of Definition \ref{dfn:equivalent} in the node $j$
if $\ga_j$ is non compact. In a finite similar steps one can
obtain $\gG_2$.

Therefore we can conclude the following.

\begin{prop} \label{prop:F}
There is a map $\CF$ from the set of $K$-conjugacy classes of Lie subalgebras generated by KS-triples to the set of equivalent classes of abstract weighted Vogan diagrams.
\end{prop}

\begin{cor} \label{cor:F}
There is a map $\CF_\kp$ from the set of nilpotent $K$-orbits of $\kp$ to the set of equivalent classes of abstract weighted Vogan diagrams. Moreover, the  composition of  $\CF_\kp$ with the Kostant-Sekiguchi correspondence, gives a map $\CF_{\kg_o}$ from the set of nilpotent $G_o$-orbits of $\kg_o$ to the set of equivalent classes of abstract weighted Vogan diagrams.
\end{cor}

\begin{proof}
Let $\CO$ be a nilpotent $K$-orbit in $\kp$. According to Theorem
\ref{thm:Noelcorresp}, there is a $K$-conjugacy class of a pair
$(\kq_\km,\eta_\km) \in \CL_\km$  corresponding to $\CO$ where
$\km$ is a $\gth$-stable Levi subalgebra of $\kg$. By Theorem
\ref{thm:noticedcorresp}, this pair  consists in the $\gth$-stable
parabolic subalgebra associated to a noticed KS-triple $(h,e,f)$
in $[\km,\km]$ such that $\CO = K.e$ and $\eta_\km= \km^{(2)}\cap
\kp$. As these correspondences are one-to-one, we can define
$\CF_\kp(\CO)$ as the image by $\CF$ of the $K$-conjugacy class of
the Lie subalgebras generated by $(h,e,f)$.
\end{proof}

\begin{dfn}
A {\it weighted Vogan diagram} is an element of a class in the image of $\CF$ or $\CF_\kp$ or  $\CF_{\kg_o}$ (all these images are the same).
\end{dfn}

\begin{example}
The first diagram of Figure \ref{fig:equivalent} of the previous
section gives place to a 4-tuple $(\ks\ko(2,5),\gth,\kh_o, \gD^+)$
and the second one to $(\ks\ko(2,5),\gth,\kh_o, s_{\ga_2}(\gD^+))$
where $s_{\ga_2}$ is the reflexion associated to the simple root
${\ga_2}$ of $\gD^+$.
\end{example}

\begin{prop} \label{prop:underlying}
Let $\gG=(D,\gth_\gG,J,\go)$ be the weighted Vogan diagram  corresponding to a nilpotent $G_o$-orbit $\CO_o$ of $\kg_o$. Then,

$(i)$ the  underlying weighted Dynking diagram $(D,w)$ of $\gG$  is the weighted Dynkin diagram of the complex nilpotent $G$-orbit $\CO=G.\CO_o$;

$(ii)$ the underlying Vogan diagram $(D,\gth_\gG,J)$ of $\gG$ is a Vogan diagram of $\kg_o$.
\end{prop}

\begin{proof}
Starting with a nilpotent $G_o$-orbit $\CO_o$ of $\kg_o$ and fixing a Cartan involution $\gth$ of $\kg_o$, we can associate to it a nilpotent $K$-orbit of $\kp$ by the Kostant-Sekiguchi correspondence (see Section \ref{data}). The real nilpotent orbit $\CO_o=G_o.e_o$ is related with the nilpotent $K$-orbit $\CO_\kp=K.e$ by a Cayley transform (see 1.4), that is $e=c(e_o)=g.e_o$ for a particular element  $g\in G$. So, $G.\CO_\kp= G.\CO_o= \CO$. Then, the weighted Dynkin diagrams associated to the KS-triple $(h,e,f)=(c(h_o),c(e_o),c(f_o))$ and the real JM-triple $(h_o,e_o,f_o)$ are the same.

On the other hand,  the weighted Vogan diagram $\gG$ is the  one
associated to $(h,e,f)$. Following the proof of the Existence
Theorem of real semisimple Lie algebras of Vogan diagrams (see
Theorem 6.88, \cite{K}), the involution   $\gth_\gG$ of $D$ give
place to an involution of $\kg$, observe that it coincides with
$\gth$ by construction of $\gth_\gG$. As we know, the compact real
form of $\kg$ is $\kk_o\oplus i\kp_o$. Then, the real semisimple
Lie algebra associated to  $(D,\gth_\gG,J)$ is $\kk_o\oplus \kp_o
=\kg_o$ as we wanted.
\end{proof}

\begin{rmk}
According to the proof of last proposition, there is no confusion
to denote both,  the involution of a weighted Vogan diagram
$\gG=(D,\gth,J,\go)$ and the involution of $\kg$, by $\gth$.
\end{rmk}

\subsection{}Given a weighted Vogan diagram denote by $j_0=0$, $j_{m+1}= n+1$ and by $j_1,\dots , j_m$ the nodes such that $\go_{j_i}\neq 0$ and $j_1 <\dots < j_m$.

\begin{prop} \label{prop:1node}
Each equivalent class of weighted Vogan diagrams contains a diagram $\gG$ with the following property:

{\rm (P)} each weighted Vogan subdiagram $\gG_{j_i,j_{i+1}}$ of
$\gG$ with nodes $ j_i +1, \dots ,j_{i+1}-1$ has at most one
painted node, for $i= 0, 1, \dots, m$.
\end{prop}

\begin{proof}
Let $\gG$ be a weighted Vogan diagram and $\gG_{j_i,j_{i+1}}$ the subdiagram defined above. The Weyl group associated to $\gG_{j_i,j_{i+1}}$ is isomorphic to a subgroup $W_{j_i,j_{i+1}}$ of $W_o=\rm{Span}\{s_\ga\in W : \gth \ga=\ga,\  \go_\ga=0\}$. By a result in \cite{K} there is a Vogan diagram with at most one painted node associated to the same real Lie algebra than the underlying Vogan diagram of $\gG_{j_i,j_{i+1}}$. The 4-tuple associated to them differ in the systems of positive roots by an element $s_i$ of $W_{j_i,j_{i+1}}$. The element $s_i$ is a composition of reflexions in $W_o$. As we have seen before Proposition \ref{prop:F}, the action on the diagram $\gG$ is nothing more than an application of finite operations of type (A).

Using this process for each $i=0,\dots, m$ we can conclude that the resulting weighted Vogan diagram has property (P).
\end{proof}

\section {Weighted Vogan diagrams of noticed nilpotent orbits}\label{diagroots}
In this section we will give a characterization of weighted Vogan
diagrams corresponding to noticed nilpotent $G_o$-orbits of
$\kg_o$.

\subsection{}
According with Subsection \ref{equivalent}, the following results are direct consequences of Lema \ref{lem:g(j)}.

\begin{lem} \label{lem:diagroots1}
Let $\gG=(D,\gth_\gG,J,\go)$ be a weighted Vogan diagram  and
$(\kg_o,\gth,\kh_o,\gD_\gG^+)$ be the 4-tuple associated to it.
Then, for each integer $j\leq n$, the following conditions are
equivalent,
\begin{enumerate}
\item[$(i)$] there exists a root $\ga=\sum_{i=1}^n n_{i}\ga_{i} \in \gD_\gG^+$ with $n_j>0$ such that $\kg_\ga \subset \kg^{(2)}\cap \kp$;

\item[$(ii)$] for the $j$-node of $\gG$ there is a connected weighted Vogan subdiagram $\gG^j=(D^j,\gth^j,J^j,\go^j)$ that contains it and satisfies one of the following conditions,
\begin{enumerate}
\item[(a)] $0<\sum_i \go^j_i \leq 2$ and there is a root $\ga \in \gD_{\gG^j}^+\subset  \gD_{\gG}^+$ of weight 2 with odd painted lenght $m_\ga$ and $n_j >0$,  or
\item[(b)] $\gG^j$ is of the type
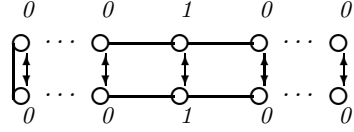
\begin{center}
\setlength{\unitlength}{10pt}
\begin{picture}(18,7)(-6,-3)
  \thicklines
  \put(-3,1){\circle{.6}}
  \put(-0,1){\circle{.6}}
  \put(3,1){\circle{.6}}
  \put(6,1){\circle{.6}}
  \put(9,1){\circle{.6}}
  \put(-3,-1){\circle{.6}}
  \put(0,-1){\circle{.6}}
  \put(3,-1){\circle{.6}}
  \put(6,-1){\circle{.6}}
  \put(9,-1){\circle{.6}}

  \put(-2.5,1){\ \dots}
  \put(0.3,1){\line(3,0){2.5}}
  \put(3.3,1){\line(3,0){2.5}}
  \put(6.5,1){\ \dots}
  \put(-3.3,1){\line(0,-1){2.1}}
  \put(-2.5,-1){\ \dots}
  \put(0.3,-1){\line(3,0){2.5}}
  \put(3.3,-1){\line(3,0){2.5}}
  \put(6.5,-1){\ \dots}

  \put(-3,2){\small{0}}
  \put(0,2){\small{0}}
  \put(3,2){\small{1}}
  \put(6,2){\small{0}}
  \put(9,2){\small{0}}
  \put(-3,-2){\small{0}}
  \put(0,-2){\small{0}}
  \put(3,-2){\small{1}}
  \put(6,-2){\small{0}}
  \put(9,-2){\small{0}}

  \thinlines
  \put(0.2,0.5){\vector(0,-1){1.1}}
  \put(0.2,-0.5){\vector(0,1){1.1}}

  \put(3.2,0.5){\vector(0,-1){1.1}}
  \put(3.2,-0.5){\vector(0,1){1.1}}

  \put(-2.8,0.5){\vector(0,-1){1.1}}
  \put(-2.8,-0.5){\vector(0,1){1.1}}

  \put(6.2,0.5){\vector(0,-1){1.1}}
  \put(6.2,-0.5){\vector(0,1){1.1}}

  \put(9.2,0.5){\vector(0,-1){1.1}}
  \put(9.2,-0.5){\vector(0,1){1.1}}

\end{picture}
\end{center}

\end{enumerate}
\end{enumerate}

\end{lem}

\begin{proof}
Suppose $\ga$ is as in $(i)$. Then, $\ga \in M_\kp^{(2)}$ by Lemma
\ref{lem:g(j)} . Let $N_\ga^j =\{{i_l}: n_{i_l} >0 \}$. Denote by
$\gG^j=(D^j,\gth^j,J^j,\go^j)$ the  connected weighted Vogan
subdiagram such that $D^j$ is the Dynkin diagram supported on
$N_\ga^j$, $J^j = J\cap N_\ga^j$, $\gth^j (i) = \gth_\gG (i)$  and
$\go^j_i = \go_i$ for all $i\in N_\ga^j$. In particular,  $\ga \in
\gD_{\gG^j}^+$. By Lemma \ref{lem:g(j)}, $m_\ga$ is odd or
$\ga=\gga+\gth\gga$ because $\kg_\ga \subset \kg^{(2)}\cap \kp$.
Then, $\ga$ satisfies $(ii.$a) or $\gG^j$ is as in $(ii.$b).

Conversely, given $\gG^j$ that satisfies $(ii.$a), consider $\ga$ a root of
$\gD_\gG^+$ as in $(ii.$a). Then, $\ga \in M_\kp^{(2)}$. It  implies that $\kg_\ga \subset \kg^{(2)}\cap \kp$ by Lemma \ref{lem:g(j)}.

If $\gG^j$ is as in $(ii.$b), denote by $\Psi^j$ the subset of simple roots associated to $D^j$ and define the root $\ga = \sum_{\ga_i \in \Psi^j} \ga_{i}$. Then, $\ga=\gga+\gth\gga$, so it is in $M_\kp^{(2)}$. Applying again Lemma \ref{lem:g(j)}, the proof is finished.
\end{proof}

\begin{lem} \label{lem:diagroots2}
Let $\gG=(D,\gth_\gG,J,\go)$ be a weighted Vogan diagram  and
$(\kg_o,\gth,\kh_o,\gD_\gG^+)$ be the 4-tuple associated to it.
Then, for each integer $j\leq n$, the following conditions are
equivalent,

$(i)$  there exists a root $\ga \in \gD_\gG^+$ with $n_j>0$ such that $\kg_\ga \subset \kg^{(2)}\cap \kk$;

$(ii)$ for the $j$-node of $\gG$ there is a connected weighted Vogan subdiagram $\gG^j=(D^j,\gth^j,J^j,\go^j)$ that contains it,
 $0<\sum_i \go^j_i \leq 2$ and there is a root $\ga \in \gD_{\gG^j}^+\subset  \gD_{\gG}^+$ of weight 2 such that its painted lenght $m_\ga$ is even  and $n_j >0$.
\end{lem}

The proof is analogous to the previous one.

\begin{lem} \label{lem:diagroots3}
Let $\gG=(D,\gth_\gG,J,\go)$ be a weighted Vogan diagram  and
$(\kg_o,\gth,\kh_o,\gD_\gG^+)$ be the 4-tuple associated to it.
Then, for each integer $j\leq n$,the following conditions are
equivalent,

$(i)$ there exists a complex root $\ga \in \gD_\gG^+$ with $n_j>0$
and weight $2$;

$(ii)$ for the $j$-node of $\gG$ there is a connected non-$\gth_\gG$-stable  weighted Dynkin subdiagram $\gO^j=(D^j,\go^j)$ that contains it such that
$\sum_i \go^j_i = 2$.
\end{lem}

\begin{proof}
Suppose $\ga$ is as in $(i)$. As it is complex, the set $N_\ga^j =\{{i_l}: n_{i_l} >0 \}$ is not $\gth$-stable. Define the connected weighted Dynkin subdiagram $\gO^j=(D^j,\go^j)$ with $D^j$ supported in $N_\ga^j$ and $\go^j_i=\go_i$ for all $i \in N_\ga^j$. Observe that $\gth_\gG$ is not the identity automorphism. Regarding diagrams with this property, the only posibilities of the numbers $n_i$ for any positive root are $0$ or $1$. Then, $\sum_{i\in N^j} \go^j_i =\go_\ga= 2$.

Conversely, denote by $\Psi^j$ the subset of simple roots associated to $\gO^j$ and define the root $\ga = \sum_{\ga_i \in \Psi^j} \ga_{i} \in \gD_\gG^+$. Then $\ga$ satisfies $(i)$ because $D$ is not $\gth$-stable and $\go_\ga= \sum_{i\in N_\ga^j} \go^j_i = 2$.
\end{proof}


\begin{prop} \label{prop:diag<g2}
Let $\gG=(D,\gth_\gG,J,\go)$ be a weighted Vogan diagram
corresponding to the KS-triple $(h,e,f)$ of $\kg$ and let
$(\kg_o,\gth,\kh_o,\gD_\gG^+)$ be the 4-tuple associated to it.
Then, the following statements are equivalent,
\begin{enumerate}
\item $\kg$ is the minimal $\gth$-stable Levi subalgebra that contains $\kh \oplus \kg^{(2)}\cap \kp$;

\item one of the next conditions is satisfied  for each node $j$ of $\gG$,
\begin{enumerate}
\item there is a connected non-$\gth$-stable weighted Dynkin subdiagram $\gO^j=(D^j,\go^j)$ that contains the node $j$ such that
$\sum_i \go^j_i = 2$, or
\item there is a connected weighted Vogan subdiagram $\gG^j=(D^j,\gth^j,J^j,\go^j)$ that contains the node $j$ such that
\begin{enumerate}
\item[I.] $0<\sum_i \go^j_i \leq 2$ and there is a root $\ga \in \gD_{\gG^j}^+\subset  \gD_{\gG}^+$ of weight 2 with odd painted lenght $m_\ga$ and $n_j >0$,  or
\item[II.] $\gG^j$ is
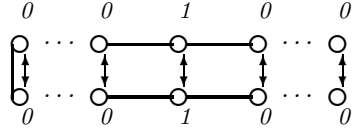
\begin{center}
\setlength{\unitlength}{10pt}
\begin{picture}(18,7)(-6,-3)
  \thicklines
  \put(-3,1){\circle{.6}}
  \put(-0,1){\circle{.6}}
  \put(3,1){\circle{.6}}
  \put(6,1){\circle{.6}}
  \put(9,1){\circle{.6}}
  \put(-3,-1){\circle{.6}}
  \put(0,-1){\circle{.6}}
  \put(3,-1){\circle{.6}}
  \put(6,-1){\circle{.6}}
  \put(9,-1){\circle{.6}}

  \put(-2.5,1){\ \dots}
  \put(0.3,1){\line(3,0){2.5}}
  \put(3.3,1){\line(3,0){2.5}}
  \put(6.5,1){\ \dots}
  \put(-3.3,1){\line(0,-1){2.1}}
  \put(-2.5,-1){\ \dots}
  \put(0.3,-1){\line(3,0){2.5}}
  \put(3.3,-1){\line(3,0){2.5}}
  \put(6.5,-1){\ \dots}

  \put(-3,2){\small{0}}
  \put(0,2){\small{0}}
  \put(3,2){\small{1}}
  \put(6,2){\small{0}}
  \put(9,2){\small{0}}
  \put(-3,-2){\small{0}}
  \put(0,-2){\small{0}}
  \put(3,-2){\small{1}}
  \put(6,-2){\small{0}}
  \put(9,-2){\small{0}}

  \thinlines
  \put(0.2,0.5){\vector(0,-1){1.1}}
  \put(0.2,-0.5){\vector(0,1){1.1}}

  \put(3.2,0.5){\vector(0,-1){1.1}}
  \put(3.2,-0.5){\vector(0,1){1.1}}

  \put(-2.8,0.5){\vector(0,-1){1.1}}
  \put(-2.8,-0.5){\vector(0,1){1.1}}

  \put(6.2,0.5){\vector(0,-1){1.1}}
  \put(6.2,-0.5){\vector(0,1){1.1}}

  \put(9.2,0.5){\vector(0,-1){1.1}}
  \put(9.2,-0.5){\vector(0,1){1.1}}

\end{picture}
\end{center}
\end{enumerate}

\end{enumerate}
\end{enumerate}
\end{prop}

\begin{proof}
Let $\km$ be a minimal $\gth$-stable Levi subalgebra of $\kg$ that
contains $\kh \oplus\kg^{(2)}\cap \kp$. Then, as $\km$ contains
$\kh$ the roots system $\gD(\km,\kh)$ is a subsystem of
$\gD_\gG(\kg,\kh)$. Hence, by Lemmas \ref{lem:diagroots3},
\ref{lem:diagroots1} and \ref{lem:g(j)},  the node $j$ satisfy
condition (a) or (b) if and only iff there is $\ga \in \gD_\gG^+$
with $n_j \neq 0$ such that $(\kg_\ga\oplus \kg_{\gth\ga}) \cap
\kg^{(2)}\cap \kp \neq \emptyset$ or $\kg_\ga \subset
\kg^{(2)}\cap \kp$ respectively. Then, this happens for each node
$j$ of $\gG$ if and only if every simple root $\ga_j$ of $\gD_\gG$
is in $\gD(\km,\kh)$, or equivalently, if and only if $\km=\kg$.
\end{proof}

Following the notation of \ref{abnoticed} we have the following
results that caracterize weighted Vogan diagrams associated to
noticed KS-triples of $\kg$.

\begin{thm} \label{thm:necesary}
If $\gG$ is a weighted Vogan diagram associated to a noticed
KS-triple $(h,e,f)$ of $\kg$, then
\begin{enumerate}
\item $P_{np}^{(1)}(\gG)$ and $P_p^{(1)}(\gG)$ have the same cardinality;
\item the statement (2) of Proposition \ref{prop:diag<g2} holds.
\end{enumerate}
\end{thm}

\begin{proof}
The statement (1) follows inmediatly from Proposition
\ref{cor:M1}, since the sets $P_{np}^{(1)}$ and $P_p^{(1)}$ are
equal to the sets $M_\kk^{(1)}$ and $M_\kp^{(1)}$ defined in
section \ref{g(i)} for the system root associated to the KS-triple
$(h,e,f)$.

 Suppose that $(h,e,f)$ is a noticed KS-triple of
$\kg$, this means that $\kg$ is the minimal $\gth$-stable Levi
subalgebra of $\kg$ that contains $e$. But $e \in \kg^{(2)}\cap
\kp$, then $\kg$ is contained in the minimal $\gth$-stable Levi
subalgebra $\km$ of $\kg$ that contains $\kh \oplus\kg^{(2)}\cap
\kp$. So, (1) of Proposition \ref{prop:diag<g2} is prooved,  or
equivalently, (2) is true. This proves the second statement.
\end{proof}

\begin{thm}
The following statements are equivalent,
\begin{enumerate}
\item  $(h,e,f)$ is a noticed KS-triple of $\kg$;

\item $\CO_\kp= K.e$ is a noticed nilpotent $K$-orbit of $\kp$;

\item $\CO_o= G_o.e_o$ is a noticed nilpotent $G_o$-orbit of $\kg_o$, where $e_o$
is the corresponding element of $e$ by the Kostant-Sekiguchi correspondence;

\item the weighted Vogan diagram associated to $(h,e,f)$ is a
noticed weighted Vogan diagram.
\end{enumerate}
\end{thm}

\begin{proof}
The three first items are equivalent by Definition \ref{dfn:
noticed}.

By Corollary \ref{cor:d0k=d2p}, $(h,e,f)$ is noticed if and only
if the subsets of the positive root system $\gD^+$ associated to
$(h,e,f)$ satisfy $\dim\kt + 2|M_\kk^{(0)}| + 2|C^{(0)}|=
|M_\kp^{(2)}| + |C^{(2)}|$. But it is obvious that they are
related to $\gG $ in the following way: $\dim\kt= N^\gth
+N^\gth_2$, $M_\kk^{(0)}=P_{np}^{(0)}(\gG)$,
$M_\kp^{(2)}=P_p^{(2)}(\gG)$, $C^{(0)}=K^{(0)}(\gG)$ and
$C^{(2)}=K^{(2)}(\gG)$. So, the equality in terms of the sets
defined from $\gG$ is exactly the condition on $\gG$ to be
noticed. So, (1) and (4) are equivalent.
\end{proof}

Given a weighted Vogan diagram, the advantage of this
caracterization is that it permits to decide easily if it is the
associated one to a noticed nilpotent orbit or not. Let's see some
examples.

\begin{example}
a) The first diagram of Figure \ref{fig:equivalent} is not a
noticed weighted Vogan diagram because the first node does not
satisfy the condition (2) of Proposition \ref{prop:diag<g2}. Other
argument is that the following numbers are are not equal,
$$
\begin{aligned}
 &N^\gth +N^\gth_2 + 2|P_{np}^{(0)}| + 2|K^{(0)}|= 3+0 + 0 + 0 \\
 &|P_p^{(2)}| + |K^{(2)}| = |\{\ga_2+2\ga_3\}| + 0=1
\end{aligned}
$$

b) The diagram of Figure \ref{fig:D6} is not a noticed weighted
Vogan diagram. In fact, beside of $\gG$ satisfies (1) and (2) of
Theorem \ref{thm:necesary}, we obtain that
$$
\begin{aligned}
 &N^\gth +N^\gth_2 + 2|P_{np}^{(0)}| + 2|K^{(0)}|= 4+1 + 6 + 0 =11\\
 &|P_p^{(2)}| + |K^{(2)}| = \\
 &\qquad=|\{\ga_1, \ga_1+\ga_2, \ga_1+\ga_2+\ga_3, \ga_4+\ga_5+\ga_6, \\
 &\qquad \qquad \ga_3+\ga_4+\ga_5+\ga_6, \ga_2+\ga_3+\ga_4+\ga_5+\ga_6\}| + 0=6
\end{aligned}
$$

c) The following diagram is a noticed weighted diagram,
\begin{figure}[h]\label{fig:nilpotent}
\begin{center}
\setlength{\unitlength}{10pt}
\begin{picture}(12,7)(-6,-3)
  \thicklines
  \put(-6,0){\circle{.6}}
  \put(-3,0){\circle*{.6}}
  \put(0,0){\circle{.6}}
  \put(3,0){\circle*{.6}}
  \put(6,0){\circle*{.6}}
  \put(-5.7,0){\line(3,0){2.5}}
  \put(-2.7,0){\line(3,0){2.5}}
  \put(.3,0){\line(3,0){2.5}}
  \put(3.2,-0.2){\line(3,0){2.6}}
  \put(3.2,0.2){\line(3,0){2.6}}
  \put(4.5,-0.2){$\rangle$}

  \put(-6,1){\small{2}}
  \put(-3,1){\small{0}}
  \put(0,1){\small{0}}
  \put(3,1){\small{2}}
  \put(6,1){\small{0}}
\end{picture}
\end{center}
\caption{} \label{nilp-fig}
\end{figure}

In fact,
$$
\begin{aligned}
 &N^\gth +N^\gth_2 + 2|P_{np}^{(0)}| + 2|K^{(0)}|= 5+0 + 2 + 0 =7\\
 &|P_p^{(2)}| + |K^{(2)}| = \\
 &\qquad=|\{\ga_1+\ga_2, \ga_1+\ga_2+\ga_3, \ga_4, \ga_4+2\ga_5, \ga_3+\ga_4, \\
 &\qquad \qquad   \ga_2+\ga_3+\ga_4+\ga_5,\ga_3+\ga_4+2\ga_5\}| +
 0=7.
\end{aligned}
$$
Moreover, this diagram is equivalent to the diagram with the same
weights but with all the nodes painted or to the one with the
first, third and fifth nodes painted. In this case there more than
one diagram with the property (P) of Proposition \ref{prop:1node}.
\end{example}

\

In a future paper we will present the weighted Vogan diagrams
associated to noticed nilpotent orbits.

\providecommand{\MR}{\relax\ifhmode\unskip\space\fi MR }
\providecommand{\MRhref}[2]{%
  \href{http://www.ams.org/mathscinet-getitem?mr=#1}{#2}
}
\providecommand{\href}[2]{#2}

\end{document}